\documentclass[11pt]{amsart}
\usepackage{epsfig}
\usepackage{amsmath}
\usepackage{amssymb}
\usepackage{amsthm}
\newcommand{\bm}[1]{\mbox{\boldmath $#1$}}

\newtheorem{theorem}{Theorem}[section]

\newtheorem{conjecture}[theorem]{Conjecture}
\theoremstyle{definition}
\newtheorem{definition}[theorem]{Definition}

\begin{document}

 \title
{On plethysm conjectures of Stanley and Foulkes: the $2 \times n$ case}

  \author {Pavlo Pylyavskyy}
\address{Department of Mathematics, Massachusetts Institute of
  Technology,
Cambridge, MA 02141}
\email{pasha@mit.edu}

\begin{abstract}
We prove Stanley's plethysm conjecture for the $2 \times n$ case, which composed with the work of Black and List provides another proof of Foulkes conjecture for the $2 \times n$ case. We also show that the way Stanley formulated his conjecture, it is false in general, and suggest an alternative formulation.
\end{abstract}

  \maketitle 

\section {Introduction}

%The wreath product of two symmetric groups $S_m w S_n = (\underbrace {S_m \times \cdots \times S_m}_{n \text{ times}}) \rtimes S_n$ is a subgroup of $S_{mn}$ in a natural way%. Consider a character $\pi_{m,n}$ of permutation action of $S_{mn}$ on cosets of $S_m w S_n$. It decomposes into sum $\sum_{\lambda} a_{\lambda}^{m,n} s_{\lambda}$ of irredu%cable characters of $S_{mn}$ indexed by partitions $\lambda$ of $mn$. A conjecture of Foulkes \cite{Foulk} says that the multiplicity $a_{\lambda}^{m,n}$ of the Schur functio%n $s_{\lambda}$ in this decomposition is less than or equal to its multiplicity $a_{\lambda}^{n,m}$ in $\pi_{n,m}$ for all partitions $\lambda$ of $nm$ when $m\geq n$.

Denote by $V$ a finite-dimensional complex vector space, and by $S^mV$ its $m$-th symmetric power. Foulkes in \cite{Foulk} conjectured that the $GL(V)$-module $S^n(S^mV)$ contains the $GL(V)$-module $S^m(S^nV)$ for $n \geq m$. For $m = 2, 3$ and $4$ the conjecture was proved; see \cite{Th}, \cite{DS}, \cite{B}. An extensive list of references can be found in \cite{V}.

In \cite{BL} Black and List showed that Foulkes conjecture follows from the following combinatorial statement. Denote $I_{m,n}$ to be the set of dissections of $\{1, \ldots, mn\}$ into sets of cardinality $m$. Let $s = \bigsqcup_{i=1}^n S_i$ and  $t = \bigsqcup_{i=1}^m T_i$ be elements of $I_{m,n}$ and $I_{n,m}$ respectively. Define matrix $M^{m,n} = (M_{t,s}^{m,n})$ by 
$$
M_{t,s}^{m,n} = 
\begin{cases}
1 & \text{if $|S_i \cap T_j|=1$ for any $1 \leq i \leq n, 1 \leq j \leq m$;}\\
0 & \text{otherwise.}
\end{cases}
$$

\begin{theorem}[Black, List 89]
If the rank of $M^{m,n}$ is equal to $|I_{n,m}|$ for integers $n \geq m > 1$, then Foulkes conjecture holds for all pairs of integers $(n,r)$ such that $1 \leq r \leq m$.
\end{theorem}

Let $\lambda$ be a partition of $N$. A {\it {tableau}} is a filling of a Young diagram of shape $\lambda$ with numbers from $1$ to $N$, and let $T_{\lambda}$ to be the set of such tableaux. Define two tableaux to be $h$-equivalent, denoted $\equiv_h$, if they can be obtained one from the other by permuting elements in rows and permuting rows of equal length. Define a {\it {horizontal tableau}} to be an element of $H_{\lambda} := T_{\lambda}/\equiv_h$. In other words, rows of a horizontal tableau form a partition of the set $\{1, \ldots, N\}$. Similarly, define $v$-equivalence $\equiv_v$ and the set $V_{\lambda} := T_{\lambda}/\equiv_v$ of {\it {vertical tableaux}} of shape $\lambda$. Consider a horizontal tableau $\mu$ with rows $r_1, \ldots, r_p$ and a vertical tableau $\nu$ with columns $c_1, \ldots, c_q$. Call $\mu$ and $\nu$ {\it {orthogonal}}, denoted $\mu \perp \nu$, if the inequality $|r_i \cap c_j| \leq 1$ holds for all $i, j$ . Equivalently, $\mu$ and $\nu$ are orthogonal if and only if there exists a tableau $\rho$ consistent with both $\mu$ and $\nu$.

Define the matrix $K_{\lambda} = (K_{\lambda}^{\mu, \nu})$ by 
$$
K_{\lambda}^{\mu, \nu} = 
\begin{cases}
1 & \text{if $\mu \perp \nu$;}\\
0 & \text{otherwise.}
\end{cases}
$$

The rows of $K_{\lambda}$ are naturally labelled by horizontal tableaux, while the columns are labelled by vertical tableaux. Let $\lambda'$ be the conjugate partition. In \cite{Stan}, Stanley formulated a conjecture, which can be equivalently stated as follows.
 
\begin{conjecture} \label{conj1}
If $\lambda \geq \lambda'$ in dominance order, i.e. $\lambda_1 + \cdots + \lambda_i \geq \lambda'_1 + \cdots + \lambda'_i$ for all $i$, then the rows of $K_{\lambda}$ are linearly independent.
\end{conjecture}

\begin{figure}
\begin{center}
\input{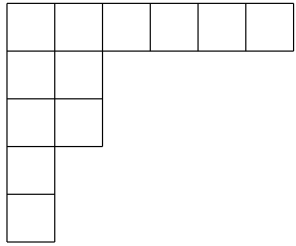}
\end{center}
\caption{A counterexample for Stanley's conjecture.}\label{lambda}
\end{figure}

This conjecture is false. For the shape $\lambda$ shown in Figure \ref{lambda}, the inequality $\lambda \geq \lambda'$ holds. However, the matrix $K_{\lambda}$ has more rows than columns, thus the rows cannot be linearly independent. Indeed, $|H_{\lambda}| = \frac{12!}{6!2!2!1!1!2!2!}$, which is greater than $|V_{\lambda}| = \frac{12!}{5!3!1!1!1!1!4!}$. This counterexample was suggested by Richard Stanley as the smallest possible one. The following conjecture seems to be a reasonable alternative formulation, although there is some evidence that it is also false in general, see \cite{N}.

\begin{conjecture} \label{conj2}
$K_{\lambda}$ has full rank for all $\lambda$.
\end{conjecture}

Let $\bm {m \times n}$ denote the rectangular shape with $m$ rows and $n$ columns. For rectangular shapes, Stanley's conjecture  implies Foulkes conjecture since $K_{\bm {m \times n}} = M^{m,n}$. For hook shaped $\lambda$, the conjecture is known to be true; see \cite{Stan}. In Section 2 we present a proof of Stanley's conjecture for $\lambda = {\bm {2 \times n}}$.

The author would like to thank Prof. Richard Stanley for the inspiration and many helpful suggestions.

\section{The Main Result}

%We begin with the following lemma.

%\begin{lemma}
%For a rectangular shape $\lambda = m \times n$ with $m$ columns and $n$ rows, $m \geq n$, we have inequality $|H_{\lambda}| \leq |V_{\lambda}|$.
%\end{lemma}

%\begin{proof}
%The statement translates to proving the inequality $\frac{(mn)!}{(m!)^{n-1}} \leq \frac{(mn)!}{(n!)^{m-1}}$, which is equivalent to $(m!)^{\frac{1}{m-1}} \geq (n!)^{\frac{1}{%n-1}}$. It is easy to check that for $m = n+1$ this holds, which implies the general case.
%\end{proof}

Our aim is to prove the following theorem.

\begin{theorem} \label{thm}
Conjecture \ref{conj2} is true for $\lambda = {\bm {2 \times n}}$.
\end{theorem}

Note that for rectangular shapes, Conjectures \ref{conj1} and \ref{conj2} are equivalent, because for $m \leq n$, the inequality $|H_{\bm {m \times n}}| \leq |V_{\bm {m \times n}}|$ holds. Therefore, proving that $K_{{\bm {2 \times n}}}$ has full rank is equivalent to proving that its rows are linearly independent. Suppose for contradiction that there is a nontrivial linear combination of rows of $K_{{\bm {2 \times n}}}$ equal to $0$. Let $\tau_{\mu}$ be the coefficient of the row corresponding to a horizontal tableau $\mu$ in this combination. Then for a column of $K_{{\bm {2 \times n}}}$ labelled by a vertical tableau $\nu$, the linear combination $\sum_{\mu} K_{{\bm {2 \times n}}}^{\mu, \nu} \tau_{\mu}$ equals $0$. Alternatively, this sum can be written as $\sum_{\mu \perp \nu} \tau_{\mu} = 0$. Call a $0$-{\it {filter}} a condition on horizontal tableax such that sum of $\tau_{\mu}$ over all $\mu$ satisfying this condition is $0$. Thus, orthogonality to $\nu$ is a $0$-filter. Our aim is to show that being $\mu$ is a $0$-filter for every horizontal tableau $\mu$. Indeed, this is just saying that all $\tau_{\mu}$ are equal to $0$, which contradicts the assumption above.

\begin{definition}
For $k<n$, a {\it {subtableau}} of shape ${\bm {2 \times k}}$ of a vertical tableau $\nu$ of shape ${\bm {2 \times n}}$ is a subset of $k$ columns of $\nu$. A {\it {partial tableau}} is a collection of $k$ columns which is a subtableau of at least one vertical tableau $\nu$.
\end{definition}

\begin{figure}
\begin{center}
\input{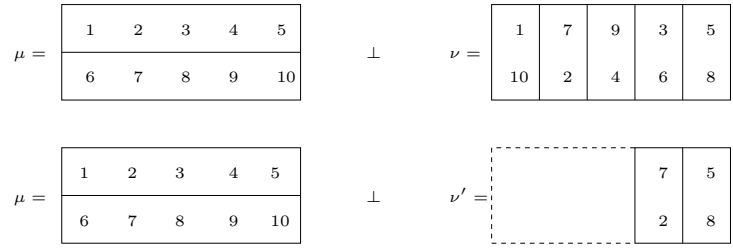}
\end{center}
\caption{Partial tableau $\nu'$ is a subtableau of $\nu$. Since $\mu \perp \nu$, also $\mu \perp \nu'$.} \label{fig2}
\end{figure}

In other words, a partial tableau is a vertical tableau of shape ${\bm {2 \times k}}$, filled with numbers from $\{1, \ldots, 2n\}$. We can now generalize the concept of orthogonality as follows. Call a horizontal tableau $\mu$ of shape ${\bm {2 \times n}}$ orthogonal to a partial tableau $\nu'$ of shape ${\bm {2 \times k}}$, where $k < n$, if there exists vertical tableau $\nu$ of shape ${\bm {2 \times n}}$ such that $\mu \perp \nu$, and $\nu'$ is a subtableau of $\nu$. An example is presented in Figure \ref{fig2}. The reason for considering such a generalization is evident from the following theorem.

\begin{theorem}
Orthogonality to a certain partial tableau $\nu'$ is a $0$-filter.
\end{theorem}

\begin{proof}
For a given partial tableau $\nu'$ of shape ${\bm {2 \times k}}$, denote $$F(\nu') = \{\nu \in V_{{\bm {2 \times n}}} \mid \text{ $\nu'$ is a subtableau of $\nu$}\}.$$ Consider the sum $\sum_{\nu \in F(\nu')} \sum_{\mu \perp \nu} \tau_{\mu}$. We claim that for each horizontal tableau $\mu$, $\tau_{\mu}$ enters this sum with the same coefficient. Indeed, the coefficient of a particular $\tau_{\mu}$ is the number of tableaux $\nu$ containing $\nu'$ and orthogonal to $\mu$. Such $\nu$'s are in one-to-one correspondance with matchings between two sets of size $n-k$, as can be seen from the Figure \ref{fig3}. The number of such matchings is $(n-k)!$, which obviously does not depend on particular $\mu$. Therefore, $\frac{1}{(n-k)!} \sum_{\nu \in F(\nu')} \sum_{\mu \perp \nu} \tau_{\mu} = \sum_{\mu \perp \nu'} \tau_{\mu}$. Since each $\sum_{\mu \perp \nu} \tau_{\mu}$ is zero by the assumption above, the sum $\sum_{\mu \perp \nu'} \tau_{\mu}$ is also $0$, which means that orthogonality to $\nu'$ is a $0$-filter.

\begin{figure}
\begin{center}
\input{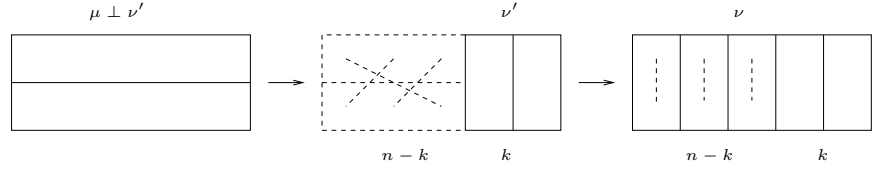}
\end{center}
\caption{Each matching corresponds to exactly one possible tableau $\nu \perp \mu$ containing $\nu'$ as a subtableau.} \label{fig3}
\end{figure}
\end{proof}

We now continue the proof of Theorem \ref{thm}. Choose a particular horizontal tableau $\mu_0$, for example with one row filled with numbers $1, \ldots, n$ and the other row filled with the rest of the numbers. If we show that $\tau_{\mu_0} = 0$, then in a similar fashion (just by relabelling numbers) we can show that all $\tau_{\mu}$'s are $0$, which would be a contradiction with the assumption that the combination of rows of  $K_{{\bm {2 \times n}}}$ is nontrivial. For a given horizontal tableau $\mu$, let $a_{\mu}$ and  $b_{\mu}$ be the numbers of elements of $\{1, \ldots, n\}$ in the rows of $\mu$, so that $a_{\mu} +  b_{\mu} = n$. We do not distinguish between the rows of $\mu$, therefore we can assume that $a_{\mu} \geq b_{\mu}$. Observe that $\mu_0$ is the only horizontal tableau such that $(a_{\mu_0}, b_{\mu_0}) = (n, 0)$. Let $T_a$ be the collection of horizontal tableaux $\mu$ with $a_{\mu} = a$, and call elements of $T_a$ {\it {horizontal tableaux of type}} $a$. Then $\mu_0$ is the only horizontal tableau of type $n$. 

\begin{theorem}
For $a \geq n/2$, being a horizontal tableau of type $a$ is a $0$-filter.
\end{theorem}

\begin{proof}
For $k \leq [n/2]$, consider the set $P_k$ of all possible partial tableaux of shape ${\bm {2 \times k}}$, filled with numbers from $\{1, \ldots, n\}$. Consider the sum $\sum_{\nu' \in P_k} \sum_{\mu \perp \nu'} \tau_{\mu}$. We claim that only $\tau_{\mu}$'s for $\mu$ of type at most $n - k$ appear in this sum. We also claim that the coefficient of $\tau_{\mu}$ in the sum depends only on the type of $\mu$.

The first statement is easy to verify. Let $\mu$ be orthogonal to some $\nu' \in P_k$. Then each of the two rows of $\mu$ contains at least $k$ numbers from $\{1, \ldots, n\}$, which means it cannot have type larger than $n-k$. As for the second statement, we can calculate exactly the number of different $\nu' \in P_k$ that are orthogonal to a given $\mu$ of type $a$. Indeed, first choose an unordered $k$-tuple among the $n-a$ elements of $\{1, \ldots, n\}$ in one row of $\mu$. Then match them with a ordered $k$-tuple taken from the $a$ elements of  $\{1, \ldots, n\}$ in the other row. Obviously, such a procedure gives all possible $\nu'$, each exactly once. Therefore, the coefficient of $\tau_{\mu}$ which we are looking for is $c_a^{k} = \frac{(n-a)!a!}{k!(n-a-k)!(a-k)!}$.

We now proceed by induction. For the base case, take $k = [n/2]$. The only horizontal tableaux in the sum $\sum_{\nu' \in P_k} \sum_{\mu \perp \nu'} \tau_{\mu}$ are those of type $n - [n/2]$. Since they all have the same coefficient, and $\sum_{\nu' \in P_k} \sum_{\mu \perp \nu'} \tau_{\mu} = 0$ because each $\sum_{\mu \perp \nu'} \tau_{\mu} = 0$, we conclude that $\sum_{\mu \in T_{n - [n/2]}} \tau_{\mu} = 0$. 

Given that being a type $a$ tableau is a $0$-filter for $n - [n/2] \leq a \leq a' < n$, let us show that being a type $a' + 1$ tableau is a $0$-filter. Indeed, $\sum_{\nu' \in P_{n-a'-1}} \sum_{\mu \perp \nu'} \tau_{\mu} = 0$ as before. This equality can be written as $\sum_{n- [n/2] \leq a \leq a'+1} c_a^{n-a'-1} \sum_{\mu \in T_a} \tau_{\mu} = 0$, where $c_a^{n-a'-1}$ is the coefficient calculated above. By the induction assumption, we know that for $n - [n/2] \leq a \leq a'$, the sum $\sum_{\mu \in T_a} \tau_{\mu}$ is $0$. Since $c_{a'+1}^{n-a'-1} \not = 0$, we conclude that $\sum_{\mu \in T_{a'+1}} \tau_{\mu} = 0$.
\end{proof}

A trivial observation to make is that for $a = n$ this theorem implies that $\tau_{\mu_0} = 0$, which leads to the desired contradiction. Therefore, rows of $K_{{\bm {2 \times n}}}$ are linearly independent, which proves Theorem \ref{thm}.

\end {document}